\documentclass[11 pt]{amsart}
\usepackage{amssymb}
\usepackage{amsmath}
\usepackage{amsfonts}
\usepackage{graphicx}
\usepackage{amsthm}
\usepackage{enumerate}
\usepackage[mathscr]{eucal}
\usepackage{verbatim}

 \theoremstyle{plain}
\newtheorem{theorem}{Theorem}[section]
\newtheorem{lemma}[theorem]{Lemma}

\numberwithin{equation}{section}
\theoremstyle{plain}

\theoremstyle{remark}

\def\bbR{{\mathbb {R}}}

\begin{document}

\date{September, 2010}

\title
{Restricting Fourier transforms of measures to curves in $\bbR^2$}

\author[]{M. Burak Erdo\u{g}an \ \ Daniel M. Oberlin}

\address
{M. B. Erdo\u{g}an \\
Department of Mathematics \\ University of Illinois \\
Urbana, IL 61801}
\email{berdogan@math.uiuc.edu}

\address
{D. M.  Oberlin \\
Department of Mathematics \\ Florida State University \\
Tallahassee, FL 32306}
\email{oberlin@math.fsu.edu}

\thanks{The first author is partially supported by the NSF grant  DMS-0900865}


\begin{abstract}
We establish estimates for restrictions to certain curves in $\bbR^2$ of the Fourier transforms 
of some fractal measures.
\end{abstract}

\maketitle

\section{Introduction}

The starting point for this note was the following observation: if
$\mu$ is a compactly supported nonnegative Borel measure 
on $\bbR^2$ which, for some $\alpha >3/2$, is $\alpha$-dimensional in the sense that 
\begin{equation}\label{main1}
\mu \big(B(y,r)\big)\lesssim r^\alpha
\end{equation}
for $y\in \bbR^2$ and $r>0$, then 
\begin{equation}\label{first}
\int_0^{\infty}|\widehat {\mu}(t,t^2 )|^2 \,dt <\infty .
\end{equation}
The proof is easy: writing $d\lambda$ for the measure given by $dt$ on the curve $(t,t^2)$, we see that
\begin{multline}\label{label}
\int_0^{\infty}|\widehat {\mu}(t,t^2 )|^2 \,dt =\iiint e^{-2\pi i (t,t^2 )\cdot (x-y)}\,
d\mu (x )\,d\mu (y )\, dt = \\
\iint \widehat{\lambda}(x -y )\,d\mu (x )\,d\mu (y ) \lesssim 
\iint |x_2 -y_2 |^{-1/2}\,d\mu (x)\,d\mu (y),
\end{multline}
where we put $x =(x_1 ,x_2 )$ if $x\in \bbR^2 $ and the inequality comes from the van der Corput estimate
$|\widehat{\lambda}(x)| \lesssim |x_2 |^{-1/2}$. 
For fixed $y$, the compact support of $\mu$ implies that 
%
\begin{equation*}
\int {|x_2 -y_2 |^{-1/2}}  \,d\mu (x ) \lesssim 
\sum_{j=0}^{\infty}2^{j/2}\mu (\{x:|x_2 -y_2 |\leq 2^{-j}\})\lesssim
\sum_{j=0}^{\infty}2^{j/2}2^{j}2^{-j\alpha}
\end{equation*}
since $\{x :|x_2 -y_2 |\leq 2^{-j}\}$ can be covered by $\lesssim 2^{j} $ balls of radius $2^{-j}$.
Clearly the last sum is finite if $\alpha > 3/2$, and then \eqref{label} is finite since $\mu$ is a finite measure. 

The simplemindedness of this argument made it seem unlikely that the index $3/2$ is best possible, and the 
search for that best index was 
the motivation for this work.
Our results here are the following theorems: 

\begin{theorem}\label{main result}
Suppose $\phi \in C^2 ([1,2])$ satisfies the estimates 
\begin{equation}\label{main1.5}
\phi ' \approx m,\,  \phi '' \approx m
\end{equation}
for some $m \geq 1$, and let $\gamma (t)=\big(t,\phi (t)\big)$. 
Suppose $\mu$ is a nonnegative and compactly supported measure on $\bbR^2$ which is $\alpha$-dimensional in the sense that 
\eqref{main1} holds. 
Then, for $\epsilon >0$, 
\begin{equation}\label{main2}
\int_1^2 |\widehat{\mu}\big(R\,\gamma (t)\big)|^2 dt 
\lesssim R^{-\alpha /2+\epsilon}\,m^{1-\alpha },
\end{equation}
when $R\ge 2$.
Here the implied constant in \eqref{main2} depends only on $\alpha$, $\epsilon$, the implied constants in \eqref{main1} and \eqref{main1.5}, and the diameter of the support of $\mu$.
\end{theorem}

\begin{theorem}\label{estimates}
Suppose $\mu$ is as in Theorem \ref{main result}, $p>1$, and 
\item{(i)} $-1 <\gamma <\alpha p-\alpha /2 -p$ if $\,1<\alpha <2$,
\item{(ii)} $-1<\gamma <-1/2$ if $\,1/2<\alpha \leq1$,
\item{(iii)} $-1<\gamma <\alpha -1$ if $\,0<\alpha\leq1/2$.

\noindent Then 
\begin{equation}\label{est1}
\int_0^{\infty}|\widehat{\mu}(t,t^p )|^2 \, t^\gamma \, dt \leq C<\infty,
\end{equation}
where $C$ depends only on $p$, the implied constant in \eqref{main1}, and the diameter of the support of $\mu$.
\end{theorem}

\begin{theorem}\label{examples}
If \eqref{est1} holds for $p>1$ and $\alpha\in (0,2)$ with $C$ as stated in Theorem \ref{estimates}, then  
\item{(i)} $-1<\gamma\leq \alpha p-\alpha /2-p$ if $\,1<\alpha <2$,
\item{(ii)} $-1<\gamma \leq -1/2 $ if $\,1/2<\alpha \leq1$, 
\item{(iii)} $-1<\gamma \leq \alpha -1$ if $\,0<\alpha \leq1/2$.
\end{theorem}

\noindent Here are some comments:

\noindent(a) Theorem \ref{main result} is a generalization of Theorem 1 in \cite{W}, which was reproved with a simpler argument in \cite{E}. As described in \S 2, the proof of Theorem \ref{main result} is just an adaptation of ideas from \cite{W} and \cite{E}.

\noindent(b) The examples which comprise the proof of Theorem \ref{examples} are similar in spirit to those
in the proof of Proposition 3.2 in \cite{W}.

\noindent(c) If $\alpha_0$ is the infimum of the $\alpha$'s for which \eqref{main1} 
implies \eqref{first} 
whenever $\mu$ is compactly supported, it follows from Theorem \ref{estimates} that 
$\alpha _0 \leq 4/3$. Then the proof of Theorem \ref{examples} and a uniform boundedness 
argument together imply that $ \alpha_0 =4/3$.

\noindent(d) Analogs of Theorem \ref{main result} have been studied for hypersurfaces in $\bbR^d$ and,
particularly, for the sphere $S^{d-1}$. See, for example, \cite{F}, \cite{M}, \cite{Sj1}, \cite{Sj2}, \cite{E}, and \cite{E2}.

The remainder of this note is organized as follows: the proof of Theorem \ref{main result} is in 
\S 2 and the proofs of Theorems \ref{estimates} and \ref{examples} are in \S3.

\section{Proof of Theorem \ref{main result}}

As mentioned above, the proof is an adaptation of ideas from \cite{W} and \cite{E}. Specifically,
with $\mu$ as in Theorem \ref{main result} and 
\begin{equation*}
\Gamma_R =\{R\,\gamma (t): 1\leq t\leq 2 \},\, 
\Gamma_{R,\delta}=\Gamma_R +B(0,R^\delta )
\end{equation*}
for $R\geq 2$ and $\delta>0$, we will modify an uncertainty principle argument from \cite{W} to show that \eqref{main2} follows from the estimate 
\begin{equation}\label{main3}
\int_{\Gamma_{R,\delta}}|\widehat{\mu}(y)|^2\,dy \lesssim R^{1-\alpha /2 +2\delta}\,m^{2-\alpha }.
\end{equation}
We will then adapt a bilinear argument from \cite{E} to prove \eqref{main3}.

So, arguing as in \cite{W}, if $\kappa\in C^{\infty}_c (\bbR^2 )$ is equal to $1$ on the support of $\mu$, then
\begin{multline}\label{main4}
\int_1^2 \big|\widehat{\mu}\big(R\,\gamma (t)\big)|^2 dt
 =\int_1^2 \Big| \int\widehat{\kappa}\big(R\,\gamma (t)-y\big)
\,\widehat{\mu}(y)\,dy\Big|^2\,dt
\lesssim \\
\int\int_1^2 \big| \widehat{\kappa}\big(R\,\gamma (t)-y\big)\big| dt\ |\widehat{\mu}(y)|^2 \,dy .
\end{multline}
If $y=(y_1 ,y_2 )$, then
\begin{multline*}
\int_1^2 \big| \widehat{\kappa}(R\,\gamma (t)-y)\big| dt \lesssim
\int_1^2 \frac{1}{\big(1+|R\,\gamma (t)-y|\big)^{10}}\, dt \lesssim \\
\frac{1}{\big(1+\text{dist} (\Gamma_R ,y)\,\big)^8}\ \int_1^2 \frac{1}{\big(1+|R\,\phi (t)-y_2 |\big)^2} \,dt.
\end{multline*}
Estimating the last integral using the hypothesized lower bound 
on $\phi'$, we see from \eqref{main4} that 
\begin{equation}\label{main4.5}
\int_1^2 |\widehat{\mu}\big(R\,\gamma (t)\big)|^2 dt \lesssim
\frac{1}{Rm}\int\frac{|\widehat{\mu}(y)|^2}{\big(1+\text{dist} (\Gamma_R ,y)\,\big)^{8}}\, dy .
\end{equation}
Now
\begin{multline*}
\int\frac{|\widehat{\mu}(y)|^2}{\big(1+\text{dist} (\Gamma_R ,y)\,\big)^8}\, dy
=\int_{\Gamma_{R,\epsilon /2}}+\sum_{j=2}^\infty \int_{\Gamma_{R,{j\epsilon /2}}\sim\Gamma_{R,{(j-1)\epsilon /2}}}
\lesssim \\
\int_{\Gamma_{R,\epsilon /2}}|\widehat{\mu}(y)|^2 \, dy
+\sum_{j=2}^\infty R^{-8(j-1)\epsilon /2}\int_{\Gamma_{R,{j\epsilon /2}}}|\widehat{\mu}(y)|^2 \, dy.
\end{multline*}
Thus \eqref{main2} follows from \eqref{main3} and \eqref{main4.5}.

Turning to the proof of \eqref{main3}, we note that by duality 
(and the fact that $\mu$ is finite)
it is enough to 
suppose that $f$, satisfying $\|f\|_2 =1$, is supported on 
$\Gamma_{R,\delta}$ and then to establish the estimate
\begin{equation}\label{main5}
\int |\widehat{f}(y)|^2\, d\mu (y) \lesssim
 \,R^{1-\alpha /2 +2\delta}\,m^{2-\alpha }.
\end{equation}
The argument we will give differs from the proof of Theorem 3 in \cite{E} only in certain technical details. But
because those details are not always obvious, and for the convenience of any reader, we will give the complete proof.

For $y\in\bbR^2$, write $y'$ for the point on the curve $\Gamma_R$ which is closest 
to $y$ (if there are multiple candidates for $y'$, choose the one with least first coordinate). Then $y' =R\gamma (t')$ for some $t'\in [1,2]$. 
For a dyadic interval $I\subset [1,2]$, define
\begin{equation*}
\Gamma_{R,\delta,I}=\{y\in \Gamma_{R,\delta}: t' \in I\},\, f_I =f\cdot\chi_{\Gamma_{R,\delta,I}}.
\end{equation*}
For dyadic intervals $I,J \subset [1,2]$, we write $I\sim J$ if $I$ and $J$ have the same length and are not adjacent but have adjacent parent intervals. The decomposition 
\begin{equation}\label{main6}
[1,2]\times [1,2]={\bigcup_{n\ge 2}} \Big(
\bigcup _{\substack{{|I|=|J|=2^{-n}}\\{I\sim J}}}(I\times J)
\Big)
\end{equation}
leads to 
\begin{equation}\label{main7}
\int |\widehat{f}(y)|^2\, d\mu (y)\leq \sum_{n\ge 2 }\ \ \sum_{\substack{{|I|=|J|=2^{-n}}\\{I\sim J}}}
\int|\widehat{f_I} (y)\widehat{f_J}(y)|\,d\mu (y).
\end{equation}
Truncating \eqref{main6} and \eqref{main7} gives 
\begin{multline}\label{main8}
\int |\widehat{f}(y)|^2\, d\mu (y)\leq \\
\sum_{4\leq 2^n \leq 
R^{1/2}}\ \ \sum_{\substack{{|I|=|J|=2^{-n}}\\{I\sim J}}}
\int|\widehat{f_I} (y)\widehat{f_J}(y)|\,d\mu (y) +\sum_{I\in \mathcal I} \int|\widehat{f_I} (y)|^2\,d\mu (y),
\end{multline}
where $\mathcal I$ is a finitely overlapping set of dyadic intervals $I$ with $|I| \approx R^{-1/2}$.

To estimate the integrals on the right hand side of \eqref{main8}, we begin 
with two geometric observations. The first of these is that if $I\subset [1,2]$ is an interval with length $\ell$, then 
\begin{equation*}
\Gamma_{R,I}\doteq \{ R\,\big( t,\phi (t)\big): t\in I\}
\end{equation*}
is contained in a rectangle $D$ with side lengths $\lesssim R\ell m ,R\ell ^2$, which we will abbreviate by saying that $D$ is a $(R\ell m) \times (R\ell ^2)$ rectangle. (To see this, note that the since the sine of the angle between vectors $(1,M)$ and $(1,M+\kappa )$ is 
\begin{equation*}
\frac{\kappa}{\sqrt{1+M^2}\sqrt{1+(M+\kappa )^2}},
\end{equation*}
it follows from \eqref{main1.5} that the angle between tangent vectors at the beginning and ending points of the curve $\Gamma_{R,I}$ is
$\lesssim \ell /m$. Since the distance between these two points is $\lesssim R \ell m$, it is clear that 
$\Gamma_{R,I}$ is contained in a rectangle $D$ of the stated dimensions.)
Secondly, we observe that if $\ell\gtrsim R^{-1/2}$, then an $R^\delta$ neighborhood of an 
$(R\ell m) \times (R\ell ^2)$ rectangle is contained in an $(R^{1+\delta}\ell m) \times (R^{1+\delta}\ell ^2)$ rectangle.
It follows that if $I$ has length $2^{-n}\gtrsim R^{-1/2}$, then the support of $f_I$ is contained in a rectangle 
$D$ with dimensions $(R^{1+\delta}2^{-n} m) \times (R^{1+\delta}2^{-2n})$.

The next lemma is part of Lemma 3.1 in \cite{E} (the hypothesis $1\leq \alpha \leq 2$ there is not necessary
for the conclusion of that lemma). To state it, we introduce some notation: $\phi$ is a nonnegative Schwartz function such that $\phi (x)=1$ for $x$ in the unit cube $Q$,  $\phi (x)=0$ if $x\notin 2Q$,
and, for each $M>0$, 
\begin{equation*}
|\widehat{\phi}|\leq C_M \sum_{j=1}^\infty 2^{-Mj}\chi_{2^j Q}.
\end{equation*}
For a rectangle $D\subset\bbR^2$, $\phi_D$ will stand for $\phi\circ b$, where $b$ is an affine mapping
which takes $D$ onto $Q$.
\begin{lemma}
Suppose that $\mu$ is a non-negative Borel measure on $\bbR^2$ satisfying 
\eqref{main1}. Suppose $D$ is a rectangle with dimensions $R_2 \times R_1$, where $R_2 \gtrsim R_1$,
and let $D_{\text{dual}}$ be the dual of $D$ centered at the origin. Then, 
if $\widetilde{\mu}(E)=\mu (-E)$,
\begin{equation}\label{main9}
(\widetilde{\mu}\ast |\widehat{\phi _D}|)(y)\lesssim R_2^{2-\alpha},\, y\in\bbR^2 
\end{equation}
and, if $K\gtrsim 1,\, y_0 \in \bbR^2$, then
\begin{equation}\label{main10}
\int_{K\cdot D_{\text{dual}}}(\widetilde{\mu}\ast |\widehat{\phi _D}|)(y_0 +y)\,dy\lesssim K^{\alpha}R_2^{1-\alpha}R_1^{-1}.
\end{equation}
\end{lemma}
\noindent Now if $I\in\mathcal I$ and $\text{supp}f_I \subset D$ as above, the identity 
$\widehat{f_I}=\widehat{f_I}\ast \widehat{\phi _D}$ implies that
\begin{equation*}
|\widehat{f_I}|\leq (|\widehat{f_I}|^2 \ast |\widehat{\phi _D}|)^{1/2}\|\widehat{\phi_D}\|_1^{1/2}\lesssim
(|\widehat{f_I}|^2 \ast |\widehat{\phi _D}|)^{1/2} 
\end{equation*}
and so 
\begin{multline}\label{main11}
\int|\widehat{f_I} (y)|^2\,d\mu (y)\lesssim
 \int  (|\widehat{f_I}|^2 \ast |\widehat{\phi _D}|)(y)\,d\mu (y)
= \\
\int |\widehat{f_I}(y)|^2 (\widetilde{\mu}\ast |\widehat{\phi_D} |)(-y)\,dy
\lesssim 
\|f_I \|_2^2 \, R^{1-\alpha /2 +2\delta}m^{2-\alpha},
\end{multline}
where the last inequality follows from \eqref{main9} and the fact that $D$ has dimensions 
$(R^{1/2 +\delta}m)\times R^\delta$ since $2^{-n}\approx R^{-1/2}$. Thus the estimate 
\begin{equation}\label{main12}
\sum_{I\in \mathcal I} \int|\widehat{f_I} (y)|^2\,d\mu (y)\lesssim 
R^{1-\alpha /2 +2\delta}m^{2-\alpha}
\sum_{I\in \mathcal I}\|f_I\|_2^2 \lesssim
R^{1-\alpha /2 +2\delta}m^{2-\alpha}
\end{equation}
follows from $\|f\|_2 =1$ and the finite overlap of the intervals $I\in\mathcal I$
(which implies finite overlap for the supports of the $f_I ,I\in\mathcal I$).

To bound the principal term of the right hand side of \eqref{main8}, fix $n$ with $4\leq 2^n \leq R^{1/2}$
and a pair $I,J$ of dyadic intervals with $|I|=|J|=2^{-n}$ and $I\sim J$. Since $I\sim J$, the support of 
$f_I \ast f_J$ is contained in a rectangle $D$ with dimensions 
$(R^{1+\delta}2^{-n}m )\times (R^{1+\delta}2^{-2n})$. For later reference, let $v$ be a unit vector in the direction of the longer side of $D$. As in \eqref{main11},
\begin{multline}\label{main13}
\int|\widehat{f_I} (y)\widehat{f_J} (y)|\,d\mu (y)\lesssim 
\int  (|\widehat{f_I}\,\widehat{f_J}| \ast |\widehat{\phi _D}|)(y)\,d\mu (y)
= \\
\int |\widehat{f_I}(y)\widehat{f_J}(y)|\, (\widetilde{\mu}\ast |\widehat{\phi_D} |)(-y)\,dy.
\end{multline}
Now tile $\bbR^2$ with rectangles $P$ having exact dimensions $C\times (C2^{-n}m^{-1}) $ 
for some large $C>0$ to be chosen later
and having shorter axis in the direction of $v$. 
Let $\psi$ be a fixed nonnegative Schwartz function satisfying $\psi (y)=1$ if $y\in Q$, 
$\widehat{\psi}(x)=0$ if $x\notin Q$, and 
\begin{equation}\label{main13.5}
\psi \leq C_M \sum_{j=1}^\infty 2^{-Mj} \chi_{2^j Q}.
\end{equation}
Since $\sum_P \psi_P^3 \approx 1$,
it follows from \eqref{main13} that if $f_{I,P}$ is defined by 
\begin{equation*}
\widehat{f_{I,P}}=\psi_P \cdot\widehat{f_I}
\end{equation*}
then 
\begin{multline}\label{main14}
\int|\widehat{f_I} (y)\widehat{f_J} (y)|\,d\mu (y)\lesssim \\
\sum_P \Big(\int |\widehat{f_{I,P}}(y)\widehat{f_{J,P}}(y)|^2\,dy\Big)^{1/2}
\Big(\int \big| (\widetilde{\mu}\ast|\widehat{\phi_D}|)(-y)\psi_P (y)\big|
^2\,dy\Big)^{1/2}.
\end{multline}

To estimate the first integral in this sum, we begin by noting that the support of $f_{I,P}$ is contained in $\text{supp}(f_I )+P_{\text{dual}}$, where $P_{\text{dual}}$ is a rectangle dual to $P$ and centered at the origin.
Let $\widetilde I$ be the interval with the same center as $I$ but lengthened by $2^{-n}/10$ 
and let $\widetilde J$ be defined similarly. 
Since $I\sim J$, it follows that $\text{dist}(\widetilde{I} ,\widetilde{J})\geq 2^{-n}/2$.
Now the support of $f_I$ is contained in $\Gamma_{R,I}+B(0,R^\delta )$
and $P_{\text{dual}}$ has dimensions $( m2^n C^{-1})\times C^{-1}$ with the longer direction at an angle 
$\lesssim2^{-n}/m$ to any of the tangents to the curve $\big(t,\phi (t)\big)$ 
for $t\in \widetilde{I}$ (or $t\in\widetilde{J}$). Recalling that $2^n \lesssim R^{1/2}$, one can
check that, if $C$ is large enough, 
\begin{equation*}
\text{supp}(f_{I,P})\subset \Gamma_{R,\widetilde{I}}+B(0,CR^\delta )
\end{equation*}
and, similarly,
\begin{equation*}
\text{supp}(f_{J,P})\subset \Gamma_{R,\widetilde{J}}+B(0,CR^\delta ).
\end{equation*}
The following lemma will be proved at the end of this section:

\begin{lemma}\label{intersection lemma} Suppose $\phi $ satisfies the estimates 
\begin{equation*}
0<\phi' \leq  m_1 \ \text{and}\  \,   \phi '' \geq m_2
\end{equation*}
with $m_1 \geq 1$ and 
\begin{equation}\label{approx}
m_1 , m_2 \approx m.
\end{equation}
Suppose that the closed intervals $\tilde I , \tilde J \subset [1,2]$ satisfy  
$\text{dist}\,(\tilde I ,\tilde J )\geq c \,2^{-n}$.
Then, for $\delta >0$ and $x\in\bbR^2$, there is the following estimate for the two-dimensional Lebesgue measure of the intersection of translates of tubular neighborhoods of $\Gamma_{R,\tilde I}$ and $\Gamma_{R, \tilde J}$:
\begin{equation}\label{intersection estimate}
\big| \, x+\Gamma_{R,\tilde I}+B(0,C R^\delta )\ \cap\  \Gamma_{R,\tilde J}+B(0,CR^\delta ) \,
\big| \lesssim {R^{2\delta} 2^n m}.
\end{equation}
The implicit constant in \eqref{intersection estimate} depends only on the 
implicit constants in \eqref{approx} and the
positive constants $c$ and $C$.
\end{lemma}

\noindent It follows from Lemma \ref{intersection lemma} that for $x\in\bbR^2$ we have
\begin{equation}\label{main15}
\big| \,x+\text{supp}(f_{I,P})\ \cap\ \text{supp}(f_{J,P})\,\big|\lesssim R^{2\delta}2^n m .
\end{equation}
Now 
\begin{equation*}
\int |\widehat{f_{I,P}}(y)\widehat{f_{J,P}}(y)|^2\,dy =
\int|\widetilde{f_{I,P}}\ast {f_{J,P}}(x)|^2\, dx
\end{equation*}
and 
\begin{multline*}
|\widetilde{f_{I,P}}\ast {f_{J,P}}(x)|\leq\int |f_{I,P}(w-x)\,f_{J,P}(w)|\,dw \leq \\
|x+\text{supp}(f_{I,P})\cap \text{supp}(f_{J,P})|^{1/2}\, \big(|\widetilde{f_{I,P}}|^2 \ast |{f_{J,P}}|^2 (x)\big)^{1/2}.
\end{multline*}
Thus, by \eqref{main15},
\begin{multline}\label{main16}
\Big(\int |\widehat{f_{I,P}}(y)\widehat{f_{J,P}}(y)|^2\,dy \Big)^{1/2}\lesssim
R^\delta 2^{n/2}m^{1/2}\Big(\int |\widetilde{f_{I,P}}|^2 \ast |{f_{J,P}}|^2 (x)\,dx \Big)^{1/2} = \\
R^\delta 2^{n/2}m^{1/2}\|f_{I,P}\|_2 \|f_{J,P}\|_2 .
\end{multline}

To estimate the second integral in the sum \eqref{main14} we use \eqref{main13.5} to observe that
\begin{equation*}
\psi_P \lesssim \sum_{j=1}^{\infty}2^{-Mj}\chi_{2^j P}.
\end{equation*}
Thus
\begin{equation*}
\int (\widetilde{\mu}\ast|\widehat{\phi_D}|)(-y)\psi_P (y) \,dy \lesssim
 \sum_{j=1}^{\infty}2^{-Mj}\int_{2^j P} (\widetilde{\mu}\ast|\widehat{\phi_D}|)(-y) \,dy .
\end{equation*}
Noting that $2^j P \subset y_P +K D_{\text{dual}}$  for some $K\lesssim R^{1+\delta}2^{-2n+j}$ 
and some $y_P \in \bbR^2$, we apply \eqref{main10} to obtain 
\begin{multline*}
\int (\widetilde{\mu}\ast|\widehat{\phi_D}|)(-y)\psi_P (y) \,dy \lesssim \\
\sum_{j=1}^{\infty}2^{-Mj}( R^{1+\delta}2^{-2n+j})^\alpha (R^{1+\delta}2^{-n}m)^{1-\alpha }(R^{1+\delta}2^{-2n})^{-1}
\lesssim 2^{-n(\alpha -1)}m^{1-\alpha}.
\end{multline*}
Since 
\begin{equation*}
(\widetilde{\mu}\ast|\widehat{\phi_D}|)(-y) \lesssim (R^{1+\delta}2^{-n}m)^{2-\alpha }
\end{equation*}
by \eqref{main9} and since $\psi_P (y)\lesssim 1$, it follows that 
\begin{equation}\label{main19}
\Big(\int \big((\widetilde{\mu}\ast|\widehat{\phi_D}|)(-y)\psi_P (y)\big)^2 \,dy \Big)^{1/2}\lesssim
R^{1-\alpha /2 +\delta (1-\alpha /2) }2^{-n/2}m^{3/2-\alpha}.
\end{equation}

Now \eqref{main16} and \eqref{main19}
imply, by \eqref{main14}, that
\begin{equation*}
\int|\widehat{f_I} (y)\widehat{f_J} (y)|\,d\mu (y)\lesssim
R^{1-\alpha /2 +\delta (2-\alpha /2)}m^{2-\alpha}\Big(\sum_P \|f_{I,P}\|_2^2 \Big)^{1/2} \Big(\sum_P \|f_{J,P}\|_2^2 \Big)^{1/2}.
\end{equation*}
Since
\begin{equation*}
\sum_P \|\widehat{f_{I,P}}\|_2^2 =
\int |\widehat{f_I}(y)|^2 \sum_P |\psi_P (y)|^2 \, dy,
\end{equation*}
it follows from $\sum_P \psi_P ^2 \lesssim 1$ that
\begin{equation*}
\int|\widehat{f_I} (y)\widehat{f_J} (y)|\,d\mu (y)\lesssim
R^{1-\alpha /2 +\delta (2-\alpha /2)}m^{2-\alpha}\|f_I \|_2 \|f_J \|_2 .
\end{equation*}
Thus 
\begin{multline}\label{main21}
\sum_{\substack{{|I|=|J|=2^{-n}}\\{I\sim J}}} \int|\widehat{f_I} (y)\widehat{f_J}(y)|\,d\mu (y)\lesssim \\
R^{1-\alpha /2 +\delta (2-\alpha /2 )}m^{2-\alpha}\sum_{\substack{{|I|=|J|=2^{-n}}\\{I\sim J}}}
\|f_I \|_2 \|f_J \|_2 \lesssim \\
R^{1-\alpha /2 +\delta (2-\alpha /2)}m^{2-\alpha}\|f\|_2^2 .
\end{multline}
Now \eqref{main5} follows from \eqref{main8}, \eqref{main12}, \eqref{main21}, and the fact that the first sum in \eqref{main8} has $\lesssim \log R$ terms.

Here is the proof of Lemma \ref{intersection lemma}:

\begin{proof} Fix $t\in \tilde I,\ s\in \tilde J$ such that 
\begin{equation}\label{intersection}
x+R\big( t,\phi (t)\big)+\overline{B(0,CR^\delta )}\ \cap\  R\big( s,\phi (s)\big)+
\overline{B(0,CR^\delta )}\not=\emptyset
\end{equation} 
and such that $t$ is minimal subject to \eqref{intersection}. 
Without loss of generality, assume that $t<s$.
Suppose that $v$ and $w$ satisfy 
\begin{equation}\label{intersection2}
x+R\big( t+w,\phi (t+w)\big)+\overline{B(0,CR^\delta )}\ \cap\  R\big( s+v,\phi (s+v)\big)+
\overline{B(0,CR^\delta )} \not=\emptyset .
\end{equation}
We will begin by observing that 
\begin{equation}\label{delta est}
w\leq \frac{8C2^n R^{\delta -1}m_1}{c\,m_2} .
\end{equation}
From \eqref{intersection} and \eqref{intersection2} it follows that 
\begin{equation}\label{intersection3}
|w-v |,\ \big|\big(\phi (s+v )-\phi (s)\big)-\big(\phi (t+w )-\phi (t)\big)\big|\leq 4CR^{\delta -1}.
\end{equation}
Now
\begin{equation}\label{intersection4}
\big(\phi (s+v )-\phi (s)\big)-\big(\phi (t+w )-\phi (t)\big)=\int_t^{t+w}\big(\phi '(u+s-t)-\phi' (u)\big)\,du
+e
\end{equation}
where the error term $e$ satisfies $|e|\leq 4CR^{\delta -1}m_1$
because of the first inequality in \eqref{intersection3} and the bound on $\phi'$.
Since $s-t\geq c2^{-n}$, the lower bound on $\phi''$ shows that the integral in \eqref{intersection4}
exceeds $ wc 2^{-n}m_2$. Thus if 
\begin{equation*}
wc2^{-n}\,m_2 > 8C R^{\delta -1}\,m_1
\end{equation*}
(that is, if \eqref{delta est} fails)
then, since $m_1\geq 1$,  \eqref{intersection4} exceeds $4CR^{\delta-1}$, contradicting \eqref{intersection3}.

To see \eqref{intersection estimate}, define $\tilde t$ by 
\begin{equation*}
\tilde t =t+\frac{8C2^n R^{\delta -1}m_1}{c\, m_2}
\end{equation*}
and note that by \eqref{delta est} the intersection in \eqref{intersection estimate} is contained in 
a translate of
\begin{equation*}
\{R\big(u,\phi (u)\big):t\leq u\leq \tilde t \}+B(0,CR^\delta )\doteq \Gamma +B(0,CR^\delta ).
\end{equation*}
Using $\phi' \lesssim m$, the length of the curve $\Gamma$ is 
$\lesssim 2^n R^{\delta } m $.
Thus $\Gamma$ is contained in $\lesssim 2^n m$ balls of radius $R^\delta$. This implies \eqref{intersection estimate}.

\end{proof}

\section{Proof of Theorems \ref{estimates} and \ref{examples}}
{\it Proof of Theorem \ref{estimates}:}
First suppose $1<\alpha <2$. 
Choose $\epsilon >0$ such that $\gamma +2\epsilon <\alpha (p-1/2)-p$. Then apply Theorem \ref{main result} with 
$\phi (t)=R^{p-1}t^p$ and $m=R^{p-1}$ to conclude that
\begin{equation*}
\int_1^2 \big|\widehat{\mu}\big(Rt, (Rt)^p \big)\big|^2 dt \lesssim R^{ -\alpha /2 +\epsilon}R^{(p-1)(1-\alpha )}
\end{equation*}
and so 
\begin{equation*}
\int_R^{2R}|\widehat{\mu}(t,t^p )|^2 \, t^\gamma \,dt\lesssim R^{-\epsilon}.
\end{equation*}
Now \eqref{est1} follows by taking $R=2^n$.

To deal with the remaining cases we note that if $d\nu$ is $dt$ on the curve $(t,R^{p-1}t^p )$, $1\leq t\leq 2$, then there is the estimate $|\widehat{\nu}(\xi )|\lesssim |\xi |^{-1/2}$. It follows from Theorem 1 in \cite{E} that 
\begin{equation*}
\int_1^2 \big|\widehat{\mu}\big(Rt, (Rt)^p \big)\big|^2 dt \lesssim R^{-\min (\alpha ,1/2)}. 
\end{equation*}
This implies the conclusions of Theorem \ref{examples} in cases (ii) and (iii) exactly as in the preceding paragraph.

{\it Proof of Theorem \ref{examples}:}
We begin by observing that
if the conclusion \eqref{est1} of  Theorem \ref{estimates} holds for $\alpha\in (0,2)$ with $C$ depending 
only on the size of the support of the nonnegative measure $\mu$ and the implied constant 
in \eqref{main1}, then the same conclusion holds (with $C$ replaced by $16\,C$) for complex measures whose
total variation measure $|\mu |$ satisfies \eqref{main1}.

We consider first the case $\alpha\in(1,2)$. Suppose $R$ is large and positive.
It is easy to check that the set 
\begin{equation*}
\big\{(t,t^p):R\leq t\leq R+\sqrt{R}\big\}
\end{equation*}
is contained in a rectangle $D$ with (approximate) dimensions $1\times R^{p-1/2}$. Let $v$ be a unit vector in the direction of the long axis of $D$ and $c_D$ be the center of $D$. Also, denote the dual of $D$ centered at the origin by 
$D_{\text{dual}}$. Note that $D_{\text{dual}}$ is a rectangle with dimensions $1\times R^{1/2-p}$ with short axis in the direction $v$.
Fix a function $\psi\in C_c^{\infty}$ supported in $D_{\text{dual}}$ such that 
$\widehat{\psi}\gtrsim R^{(p-1/2)(1-\alpha )}$ on $D$ and $\|\psi\|_{\infty}\lesssim R^{(p-1/2)(2-\alpha )}$. 
Let $T\approx R^{(p-1/2)(\alpha -1)}$ be a natural number and define $\mu$ by
\begin{equation}\label{mudef}
\mu(y)\,\dot= \, e^{2\pi i y\cdot c_D} \sum_{k=1}^T \psi (y-{k}{T^{-1}}v).
\end{equation}
It is easy to check that $|\mu |$ satisfies \eqref{main1} independently of $R$. 
Also note that 
\begin{equation*}
|\widehat{\mu}(x)|\gtrsim
R^{(p-1/2)(1-\alpha)} \, {\chi}_{D}(x) \,\big|\sum_{k=1}^T e^{-2\pi i \frac{k}{T}v\cdot (x-c_D)}\big|.
\end{equation*}
Now if 
\begin{equation*}
|\frac{1}{T}v\cdot (x-c_D)|\leq 1/4\  (\text{mod} \,1 ),
\end{equation*}
then we have
\begin{equation*}
\big|\sum_{k=1}^T e^{-2\pi i \frac{k}{T}v\cdot (x-c_D)}\big|\gtrsim T. 
\end{equation*}
Therefore there are $ N\approx R^{p-1/2}/T\approx R^{(p-1/2)(2-\alpha )}$ subrectangles $P_1,...,P_N$ of $D$ with dimensions $1\times 1/4$ whose centers are in an arithmetic progression with distance $T$ between the adjacent points such that
\begin{equation*}
|\widehat{\mu}(x)|\gtrsim 
R^{(p-1/2)(1-\alpha)} \,T \sum_{k=1}^N \chi_{P_k}(x)\approx\sum_{k=1}^N \chi_{P_k}(x).
\end{equation*}
Using this we obtain
\begin{align*}
 \int_R^{R+\sqrt{R}} |\widehat{\mu}(t,t^p)|^2 \,t^\gamma \, dt &\gtrsim R^\gamma \int_R^{R+\sqrt{R}}\sum_{k=1}^N \chi_{P_k}(t,t^p) \, dt\\
&\gtrsim R^\gamma\frac{N}{R^{p-1}}\approx R^{\gamma-\alpha p+\alpha/2+p}.
\end{align*}
This implies that $\gamma\leq \alpha p-\alpha/2-p$ and so gives the conclusion (i) of 
Theorem \ref{examples}.

The conclusion (ii) of Theorem \ref{examples} also follows from the examples just constructed: 
since the support of $\mu$ above is contained in
a ball of radius $\approx 1$, if $|\mu |$ satisfies \eqref{main1} for some $\alpha >1$, then the same is certainly true for all $\alpha\in (0,1]$. Taking $\alpha=1+\delta $ for arbitrary $\delta >0$ gives $\gamma \leq -1/2$.

To conclude, suppose $\alpha\in (0,1/2)$ and $R>0$ is large. Let $D$ be a rectangle with dimensions $R\times R^p$
which contains 
\begin{equation*}
\big\{(t,t^p):R\leq t\leq 2R\big\},
\end{equation*}
and let $v$, $C_D$, and $D_{\text{dual}}$ be as above.
Note that now $D_{\text{dual}}$ is a rectangle with dimensions $R^{-1}\times R^{-p}$ with short axis in the direction $v$.
Fix a function $\psi\in C_c^{\infty}$ supported in $D_{\text{dual}}$ and satisfying 
$\widehat{\psi}\gtrsim R^{-\alpha}$ on $D$ and $\|\psi\|_{\infty}\lesssim R^{p+1-\alpha}$. 
Fix a natural number $T$  
with $T\approx R^\alpha$ and again define $\mu$ by \eqref{mudef}. 
%
%
As before, $|\mu |$ satisfies \eqref{main1} independently of $R$ and
there are $N\approx R^p /T\approx R^{p-\alpha}$ disjoint
subrectangles $P_1,...,P_N$ of $D$ of dimensions $1\times 1/4$ 
such that
\begin{equation*}
|\widehat{\mu}(x)|\gtrsim R^{-\alpha} \, T \sum_{k=1}^N \chi_{P_k}(x)\approx\sum_{k=1}^N \chi_{P_k}(x).
\end{equation*}
As above, that leads to
\begin{align*}
 \int_R^{2R} |\widehat{\mu}(t,t^p)|^2 \, t^\gamma\, dt 
&\gtrsim R^\gamma \int_R^{2R}\sum_{k=1}^N \chi_{P_k}(t,t^p) \,dt\\
&\gtrsim R^\gamma\frac{N}{R^{p-1}}\approx R^{\gamma+p-\alpha -(p-1)}.
\end{align*}
This gives the conclusion (iii) of Theorem \ref{examples}.

\end{document}